\documentstyle[12pt]{article}

\title{Compact endomorphisms of $H^{\infty}(D)$}

\author{Joel F. Feinstein and Herbert Kamowitz}

\begin{document}
\newcommand{\dis}{\displaystyle}
\maketitle

     Let $D$ be the open unit disc and, as usual, let
$H^{\infty}(D)$ be the algebra of bounded analytic functions on $D$
with $\|f\|=\sup_{z \in D} |f(z)|.$
With pointwise addition and multiplication, $H^{\infty}(D)$ is a well
known uniform algebra. In this note we characterize the compact endomorphisms
of  $H^{\infty}(D)$ and determine their spectra.

     We show that although not every endomorphism $T$ of $H^{\infty}(D)$
has the form $T(f)(z)=f(\phi(z))$ for some analytic $\phi$ mapping
$D$ into itself, if $T$ is compact, there is an analytic function $\psi:D \rightarrow D$
associated with  $T$. In the case where $T$ is compact,  the derivative of $\psi$ at its fixed point  determines the spectrum of $T.$

The structure of the maximal ideal space
$M_{H^{\infty}}$ is well known. Evaluation at a point $z \in D$ gives rise to an element in $M_{H^{\infty}}$ in the natural way. The remainder of $M_{H^{\infty}}$ consists of singleton Gleason parts and Gleason parts which are
analytic discs. An analytic disc, $P(m)$, containing a point $m\in M_{H^{\infty}}$,
is a subset of $M_{H^{\infty}}$ for which there exists a continuous  bijection
$L_m:D \rightarrow P(m)$ such that $L_m(0)=m$ and $\hat {f}(L_m(z))$
is analytic on $D$ for each $f \in H^{\infty}(D)$.
Moreover, the map $L_m$ has the form $\dis L_m(z)=w^*\lim \frac{z+z_{\alpha}}
{1 + \overline{z_{\alpha}}z}$ for some net ${z_{\alpha}} \rightarrow m$ in the 
w*topology, whence
   $\dis \hat{f}(L_m(z))=\lim f(\frac{z+z_{\alpha}}{1+\overline{z_{\alpha}}z})$
for all $f \in H^{\infty}(D)$.
A fiber $M_{\lambda}$ over some $\lambda \in \overline{D}\setminus{D}$, is the zero
set in $M_{H^{\infty}}$ of the function $z-\lambda.$ Each part, distinct
from $D$, is contained in exactly one fiber $M_{\lambda}$. With no 
loss of generality we let $\lambda = 1.$ We recall, too, that
two elements $n_1$ and $n_2$ are in the same part if, and only if,
$\|n_1 - n_2\| < 2$, where $\| . \|$ is the norm in the
dual space $H^{\infty}(D)^*.$

    Now let $T$ be an endomorphism of $H^{\infty}(D)$, i.e. $T$ is
a (necessarily) bounded linear map of $H^{\infty}(D)$ to itself with
$T(fg)=T(f)T(g)$ for all $f, g \in H^{\infty}(D)$.  For a given $T$,
either $T$ has the form $Tf(z)=f(\omega(z))$ for some analytic map 
$\omega:D \rightarrow D$, or $Tf=\hat{f}(n)1$ for some
$n \in M_{H^{\infty}},$ or
there exists an $m \in M_{H^{\infty}}$, a net $z_{\alpha} \rightarrow m$
in the w* topology and an analytic function $\tau:D \rightarrow D$,
with $\tau(0)=0$ for which $Tf(z)=\hat{f}(L_m(\tau(z)))$ \cite{gar}.
Further, on general principles, if $T$ is an endomorphism of $H^{\infty}(D)$
there exists a w* continuous map $\phi:M_{H^{\infty}} \rightarrow 
M_{H^{\infty}}$ with $\widehat{Tf}(n)=\hat{f}(\phi(n))$ for all $n \in M_{H^{\infty}}$. In the last case,
$\phi(z)=L_m(\tau(z))$  for  $z \in D$.

   For a given endomorphism $T$, if the induced map $\phi$ maps $D$ to
itself, then $T$ is commonly called a {\em composition operator}. Compact
composition operators on $H^{\infty}$ were completely characterized
in \cite{sw}. However, in general, $L_m(\tau(z))$ need not be in $D$, 
and so not every endomorphism of $H^{\infty}(D)$ is a composition
operator.  It is these endomorphisms that we discuss here.
Trivially, for any $n \in M_{H^{\infty}} \setminus D$, the map $T$ defined by   $T:Tf(z)=\hat{f}(n)1$
is a compact endomorphism of $H^{\infty}(D)$ which is not a composition
operator.

    Now let $P(m)$ be an analytic part and let
$T$ be an endomorphism defined by $Tf(z)=\hat{f}(L_m(\tau(z)))$
as discussed above. Also suppose that $\dis \phi:M_{H^{\infty}} \rightarrow
M_{H^{\infty}}$ is such that $\widehat{Tf}=\hat{f} \circ \phi.$
Assume that $T$ is compact. We claim that
$\overline{\tau(D)}$ is a compact subset of $D$ in the Euclidean topology.
Indeed, if we regard the endomorphism $T$ as an operator
from $H^{\infty}(D)$ into $C(M_{H^{\infty}})$, then $T$ is compact
if, and only if, $\phi$ is w* to norm continuous on $M_{H^{\infty}}$
\cite{ds}.
Since $M_{H^{\infty}}$ is itself compact and connected (in the w* topology),
$\phi(M_{H^{\infty}})$ must be compact and connected in the norm
topology on $M_{H^{\infty}}$ and so $\phi$ maps $M_{H^{\infty}}$ into
a norm compact connected subset of $P(m).$
Therefore the range, $\phi(D)=L_m(\tau(D))$
is contained in a norm compact subset of $P(m)$, and further since
$L_m^{-1}$ is an isometry in the Gleason norms on $P(m)$ and $D$ \cite{jff},
$\tau(D)=L_m^{-1}(\phi(D))$ is contained in a compact subset of $D$ in the norm topology on $D$. Since
the norm, Euclidean and w* topologies on $D$ coincide, $\overline{\tau(D)}$ is a compact
subset of $D$ in these three topologies. As a consequence, $\hat{\tau}(M_{H^{\infty}})\subset D.$

   Next consider two maps of $H^{\infty}(D)$ into itself.
The first, $C_{L_m}$ is defined by $C_{L_m}(f)(z)=\hat{f}(L_m(z))$,
and the second $C_{\tau}$ by $C_{\tau}(f)(z)=f(\tau(z)).$
Then $(C_{L_m} \circ C_{\tau})(f)(z)=C_{L_m}(f \circ \tau)(z)=\widehat{f \circ \tau} (L_m(z))$ and $(C_{\tau} \circ C_{L_m})(f)(z)=\hat{f}(L_m(\tau(z)))=Tf(z).$ But if $B$ is a Banach space and $S_1$ and $S_2$ are any 
two  bounded linear maps from $B \rightarrow B$,
the spectrum $\sigma(S_1S_2)\setminus\{0\}=\sigma(S_2S_1)\setminus\{0\}$. Thus 
we see that
 $\sigma(T)\setminus\{0\}=  \sigma(C_{L_m} \circ C_{\tau})\setminus\{0\}.$

   Since $f$ is analytic on a neighborhood of range $\hat{\tau}$ which is a   
subset of $D$, a  standard functional calculus argument gives 
$\widehat{f \circ \tau }(L_m(z))=f(\hat{\tau}(L_m(z))$. If we let
$\psi(z)=\hat{\tau}(L_m(z))$ we see that $C_{L_m} \circ C_{\tau}$
is a compact composition operator in the usual sense, and  so if $z_0 \in D$ is the
unique fixed point of $\psi$, and ${\bf N}$ is the set of positive integers, 
then $\sigma(T)=\{(\psi'(z_0))^n:n \in {\bf N}\} \cup \{0,1\}.$

     To summarize, we have shown the following.

Theorem:\ If $T$ is a compact endomorphism
of $H^{\infty}(D)$, then either $T$ has one dimensional range,
i.e. $Tf=\hat{f}(n)1$ for some $n \in M_{H^{\infty}},$
or $T$ is a composition operator in the usual sense, or $T$ has the
form $Tf(z)=\hat{f}(L_m(\tau(z)))$ where  $\tau$ is described above.
In the last case, there is a compact composition operator $C_{\psi}$,
such that 
$\sigma(T)=\sigma(C_{\psi})=\{(\psi'(z_0))^n:n \in {\bf N}\} \cup \{0,1\}$ 
where $z_0 \in D$ is the unique fixed point of
$\psi.$

     We conclude with two examples showing differences between 
composition operators and general endomorphisms .

(a) With the same terminology and symbols, suppose $\hat{\tau}$ is constant 
on $P(m)$, i.e. $\hat{\tau}(P(m))=\{\hat{\tau}(m)\}.$ Since $T$
is compact, $\hat{\tau}(m)\in D.$
Then using $C_{\tau}$ and $C_{L_m}$ as before, we show
that $T^2f=\hat{f}(n)1$ for some $n \in P(m).$  Indeed, 
$(C_{L_m} \circ C_{\tau})f=f(t_0)1$ where $t_0=\hat{\tau}(m) \in D$.
Then we   see that \[T^2 f=[(C_{\tau} \circ C_{L_m}) \circ
(C_{\tau} \circ C_{L_m})]f=\]
\[[C_{\tau} \circ (C_{L_m} \circ C_{\tau}) \circ C_{L_m}]f=
[C_{\tau} \circ (C_{L_m} \circ C_{\tau})] (\hat{f} \circ L_m)=
C_{\tau}(\hat{f} (L_m(t_0)) 1)=\hat{f}(L_m(t_0)) 1.\] Letting
$n=L_m(t_0)$ gives the result.

   One way to have $\hat{\tau}$ constant on $P(m)$ is for $\tau$ to
be continuous at $1$ in the  usual sense. 

A more interesting
example, perhaps, is to define $\tau$ by $\dis \tau(z)=
\frac{1}{2} z e^{\frac{z+1}{z-1}}$, and  $m \in M_{H^{\infty}}$
as a w* limit of a real net $x_{\alpha}$
approaching $1$. Then $\dis \hat{\tau}(L_m(z))=\lim_{\alpha}  
\tau(\frac{z+x_{\alpha}}{1+\overline{x_{\alpha}}z})=0$, and so
$T^2 f=\hat{f}(m)1$ for all $f\in H^{\infty}(D)$.
     In both cases, $\sigma(T)=\{0,1\}$.

    (b) Finally,  let $\{z_n\}$ be an interpolating Blaschke sequence
approaching $1$,  $z_1=0$, with $m$ in the w* closure of $\{z_n\}$
and $B$ the corresponding Blaschke
product. If $\dis \tau(z)=\frac{1}{2}B(z)$, then it  is
well known \cite{gar} that 
$\dis (\hat{\tau} \circ L_m)'(0)=\frac{1}{2}(\hat{B} \circ L_m)'(0) \neq 0.$
This, then, is  an example of a compact endomorphism of $H^{\infty}(D)$
which is not a composition operator but whose spectrum properly contains $\{0,1\}.$

\vspace{.3in}

{\sf  School of Mathematical Sciences

 University of Nottingham 

 Nottingham NG7 2RD, England

 email: Joel.Feinstein@nottingham.ac.uk

and 

 Department of Mathematics

 University of Massachusetts at Boston 

 100 Morrissey Boulevard 

 Boston, MA 02125-3393

 email: hkamo@cs.umb.edu}

\vspace{.5in}

{\sf This research was supported by EPSRC grant GR/M31132}

\end{document}